\documentclass[11pt]{amsart}
\usepackage{url}
\usepackage{graphicx}
\usepackage{latexsym}
\usepackage{amsfonts,amsmath,amssymb}
\newtheorem{theorem}{Theorem}[section]																
\newtheorem{lemma}[theorem]{Lemma}

\usepackage{color}
\def\si{\par\smallskip\noindent}

\def\eps{\varepsilon}

\def\a{\alpha}
\def\be{\beta}

\def\part{\partial}

\def\b1{\bold 1}

\newcommand{\beq}{\begin{equation}}
\newcommand{\eeq}{\end{equation}}

\newtheorem{Theorem}{Theorem}[section]

\theoremstyle{remark}

\numberwithin{equation}{section}
\date{\today}

\begin{document}

\title{The likely maximum size of twin subtrees in a large random tree}

\author {Mikl\'os B\'ona, Ovidiu Costin, and Boris Pittel}
\address{Department of Mathematics, University of Florida, Gainesville, FL 32611}
\email{bona@ufl.edu}
\address{Department of Mathematics, Ohio State, University, Columbus, OH 43210}
\email{costin.9@osu.edu}
\address{Department of Mathematics, Ohio State, University, Columbus, OH 43210}
\email{pittel.1@osu.edu}

\maketitle

\begin{abstract} We call a pair of vertex-disjoint, induced subtrees of a rooted trees twins if they have the same counts of vertices by out-degrees. The likely maximum size of twins in a uniformly random, rooted Cayley tree of size $n\to\infty$ is studied. It is shown that the expected number of twins of size $(2+\delta)\sqrt{\log n\cdot\log\log n}$ approaches zero,
while the expected number of twins of size $(2-\delta)\sqrt{\log n\cdot\log\log n}$ approaches infinity.
\end{abstract}

\section{Introduction and the main result.}
For a given combinatorial structure $S$, there is intrinsic interest in finding the largest substructure of a given kind that is contained in $S$, or a set of two or more identical substructures
of maximum size that are contained in $S$, or the size of the largest common substructure in two random structures $S_1$ and $S_2$. For instance, finding the length of a longest increasing subsequence of a random permutation of length is the famous Ulam problem (Section 6.5 of \cite{combperm} is a basic survey of this subject). 
An early collection of problems featuring the largest common substructures in two combinatorial objects can be found in  \cite{aldous}. In a pioneering paper \cite{bryant}, Bryant, McKenzie, and Steel  proved  that the maximum size of a common subtree for two independent copies of a uniformly random binary tree with $n$ leaves is likely to be of order
$O(n^{1/2})$. Interestingly, the proof is based on a classic result that the length of the increasing subsequence of a uniformly random permutation is likely to be of order $O(n^{1/2})$. And a recent result  of Aldous 
\cite{aldous-2022}  is 
that the maximum size of this common subtree  largest common subtree (maximum agreement subtree) of two independent uniform random binary trees on $n$ leaves is likely to be of order at least $n^\beta$, where $\beta =(\sqrt{3}-1)/2\approx 0.366$. %The best known upper bound \cite{bryant} for the size of this subtree is  close to $\sqrt{n}$, and interestingly, this model is conceptually close to a uniformly random permutation. 
Pittel \cite{pittel} proved the bound $O(n^{1/2})$ for the maximum size of a subtree common to two independent copies of a random rooted tree, namely the terminal tree of a critical Galton-Watson branching process
 conditioned on the total number of leaves being $n$. 

In \cite{bona-flajolet}, B\'ona and Flajolet used a bivariate generating function approach to compute the probability that two randomly selected phylogenetic trees are isomorphic as unlabeled
trees. The requested probability was found to be asymptotic to a decreasing exponential modulated by a 
polynomial factor.  This line of research involved the study of the expected number of symmetries in such trees, which was initiated by McKeon in \cite{McKeon91}; in \cite{bona-flajolet}, 
this number in large phylogenetic trees was found to obey a limiting Gaussian distribution.

Motivated by the cited work, we study the likely maximum size of two, almost identical, rooted (fringe) subtrees of a random rooted tree, rather than a pair of random trees. (A subtree rooted at a vertex $v$ is defined as a subtree induced by $v$ and the set of all the {\it descendants\/} of $v$.) Formally, we call two rooted trees almost identical, twins, if they have the same counts of vertices by their out-degrees. The twin subtrees are automatically vertex-disjoint. In this paper the tree in question is the uniformly random (Cayley) rooted tree  on
vertex set $[n]$. Our main result is
\begin{Theorem}\label{thmStar} Let $m_n(k)$ denote the expected number of twin pairs of subtrees, with $k$ vertices each. {\bf (a)\/}
If $\delta>0$ and 
\[
K_n=\big\lfloor\exp\bigl((2+\delta)\sqrt{\log n \cdot\log(\log n)}\bigr)\big\rfloor,
\]
then  $\sum_{k\ge K_n}m_n(k)\to 0$. Consequently, with high probability, the largest size of a pair of twins is below $K_n$. {\bf (b)\/} If $\delta\in (0,2)$ and 
\[
k_n=\big\lfloor\exp\bigl((2-\delta)\sqrt{\log n \cdot\log(\log n)}\bigr)\big\rfloor,
\]
then  $m_n(k_n)\to\infty$. %Consequently, with high probability, the largest size of a pair of twins is below $K_n$
\end{Theorem}
{\bf Note.\/} {\bf (i)\/} In fact, our proof of {\bf (b)\/} consists of showing that even a {\it smaller\/} expected number of twins of maximum degree $\lfloor 2\log k_n/\log\log k_n \rfloor$ grows indefinitely as $n\to\infty$. {\bf (ii)\/} Going out on a limb,
we conjecture that the logarithm of the maximum twin size scaled by $\sqrt{\log n \cdot\log(\log n)}$ converges, in probability, to $2$.

In conclusion, a reviewer indicated to us that the our problem  is reminiscent of the study of distinct rooted subtrees
in Flajolet, Sipala, and Steyaert \cite{Fla} and the study of repeated subtrees in Ralaivaosaona and Wagner \cite{Ral}.
\section{The proof of Theorem \ref{thmStar}.}
%{\bf Definition.\/} We call two rooted subtrees of the random rooted Cayley tree twins if they have the same counts of vertices of each out-degree. 
Given $k$, $\bold r=(r_0,r_1,\dots)$ is a sequence of counts of vertices by out-degree of a rooted tree on $[k]$ if and only if
\begin{equation}\label{1}
\sum_i r_i=k,\quad \sum_i ir_i=k-1,
\end{equation}
and $M(\bold r)$, the number of such trees, is given by
\begin{equation}\label{????}
M(\bold r)=\frac{(k-1)!}{\prod\limits_{j \ge 0}(j!)^{r_j}}\binom{k}{r_0,r_1,\dots},
\end{equation}
Stanley \cite{stanley}. It follows that the total number of pairs of rooted subtrees, on two respective vertex-disjoint vertex sets, each of cardinality $k$, 
that have the same counts of vertices by out-degrees, is
\begin{equation}\label{3}
N(k)=\sum_{\bold r\text{ meets }\eqref{1}}M^2(\bold r).
\end{equation}
{\bf Note.\/} If we consider only subtrees of maximum out-degree (strictly) below $d\,(\ge 2)$ then the index $i$ in the condition \eqref{1} ranges from $0$ to $d-1$.

Let $S_n(k)$ stand for the total number of pairs of vertex-disjoint rooted subtrees, each with $k$ vertices, in all the rooted trees on $[n]$ that share the counts of vertices by out-degree. Then
\begin{equation}\label{4}
S_n(k)=\binom{n}{k}\binom{n-k}{k} N(k) (n-2k)^{n-2k-1}(n-2k)^2.
\end{equation}
And the expected number of the pairs of twins of size $k$ each is given by $m_n(k)=\tfrac{S_n(k)}{n^{n-1}}$. 

The challenge is to sharply estimate $N(k)$ given by \eqref{3}. Using \eqref{????} for $M(\bold r)$,  
and a generating function $H(z):=\sum_{r\ge 0}\tfrac{z^r}{(r!)^2}$, we evaluate
\begin{multline}\label{3.1}
N(k)=\sum_{\bold r\text{ meets }\eqref{1}}M^2(\bold r)=\bigl(k! (k-1)!\bigr)^2\sum_{\bold r\text{ meets }\eqref{1}}\,\biggl(\prod_{j\ge 0} \tfrac{1}{(j!)^{r_j} r_j!}\biggr)^2\\
=\bigl(k! (k-1)!\bigr)^2 [x_1^k x_2^{k-1}]\sum_{\bold r\ge \bold 0}x_1^{\sum_{\a}r_{\a} }\,x_2^{\sum_{\be}\be\cdot r_{\be}}
\biggl(\prod_{j\ge 0} \tfrac{1}{(j!)^{r_j} r_j!}\biggr)^2\\
=\bigl(k! (k-1)!\bigr)^2 [x_1^k x_2^{k-1}]\sum_{\bold r\ge \bold 0} \prod_{j\ge 0}\tfrac{(x_1x_2^j)^{r_j}}{\bigl((j!)^{r_j} r_j!\bigr)^2}\\
=\bigl(k! (k-1)!\bigr)^2 [x_1^k x_2^{k-1}]\prod_{j\ge 0}\sum_{r_j\ge 0}\tfrac{(x_1x_2^j)^{r_j}}{\bigl((j!)^{r_j} r_j!\bigr)^2}\\
=\bigl(k! (k-1)!\bigr)^2 [x_1^k x_2^{k-1}]\prod_{j\ge 0} H\bigl(\tfrac{x_1x_2^j}{(j!)^2}\bigr).
\end{multline}
If we count only the twins of maximum out-degree below $d$, then in the above product $j$ ranges from $0$ to $d-1$. 

What can we get from \eqref{3.1} with a bare-minimum analytical work? One possibility is to use a Chernoff-type inequality
implied by \eqref{3.1}.  Note that
for a power series $p(x,y)=\sum_{i,j}a_{i,j}x^iy^j$ with nonnegative real coefficients, the inequality 
$a_{i,j} x^i y^j \leq p(x,y)$ holds for all $(i,j)$ and all $(x,y)$ where $x>0$ and $y> 0$.  Dividing by $x^iy^j$ and taking the infimum on the right-hand side, we get that
\[a_{i,j} \leq \inf_{x>0,y>0}  \frac{p(x,y)}{x^iy^j}.\]
Applying this inequality for the power series defined by the last line of \eqref{3.1}, we get

\begin{equation}\label{3.2}
N(k)\le\inf_{x_1>0, x_2>0}\tfrac{(k!(k-1)!)^2}{x_1^kx_2^{k-1}}\prod_{j\ge 0} H\bigl(\tfrac{x_1 x_2^j}{(j!)^2}\bigr).
\end{equation}
Observe that, by an elementary inequality $\binom{2r}{r}\le 2^{2r}$, for $x>0$
\begin{equation}\label{3.21}
H(x)\le \sum_{r\ge 0}\tfrac{x^r 2^{2r}}{(2r)!}\le \sum_{\ell\ge 0}\tfrac{(2x^{1/2})^{\ell}}{\ell!}=\exp(2x^{1/2}).\\
\end{equation}
%In fact, since $H(z^2)$ is a zero-order modified  Bessel function,
%\begin{equation}\label{3.22}
%H(x)=\Theta\bigl(\tfrac{\exp(2x^{1/2})}{x^{1/4}}\bigr),\quad x\to\infty.
%\end{equation}
So, \eqref{3.2} yields
\begin{multline*}
N(k)\le\inf_{x_1>0, x_2>0}\tfrac{(k!(k-1)!)^2}{x_1^kx_2^{k-1}}\exp\biggl(2\sum_{j\ge 0}\tfrac{(x_1{x_2^j})^{1/2}}{j!}\biggr)\\
=\inf_{x_1>0, x_2>0}\tfrac{(k!(k-1)!)^2}{x_1^kx_2^{k-1}}\exp\bigl(2x_1^{1/2}e^{x_2^{1/2}}\bigr).
\end{multline*}
A straightforward computation shows that the $k$-dependent factor in the RHS function attains its minimum at $x_1=k^2\exp\bigl(-\tfrac{2(k-1)}{k}\bigr)$,
$x_2=\bigl(\tfrac{k-1}{k}\bigr)^2$, and consequently
\begin{equation}\label{3.3}
N(k)=O\bigl(\tfrac{(k!)^4e^{4k}}{k^{2(k+1)}}\bigr).
\end{equation}
The last estimate is too crude, but the above values of $x_1$ and $x_2$ will be used when we switch from Chernoff's bound to a bivariate Cauchy integral, taken on the Cartesian product of two cycles, i. e. counterclock-wise oriented boundaries of two circles, of radii $x_1$ and $x_2$ respectively.
Essentially we will use a bivariate version of the saddle point method, with the above discussion allowing us to guess that $(x_1,x_2)$ is a promising approximation of the saddle point itself.

Here is the switch. By \eqref{3.1} and the Cauchy formula we have
\begin{equation}\label{3.105}
N(k)=\tfrac{\bigl(k! (k-1)!\bigr)^2}{(2\pi i)^2}\oint\limits_{C_1\times C_2}\prod_{j\ge 0}\tfrac{1}{z_1^{k+1} z_2^k}H\Bigl(\tfrac{z_1z_2^j}{(j!)^2}\Bigr)\,dz_1 dz_2,
\end{equation}
$C_s=\bigl\{z_s\in \Bbb C: z_s=x_se^{i\xi_s},\,\xi_s\in (-\pi,\pi)\bigr\}$, $s=1,2$. It follows that for all $x_s>0$
\begin{equation}\label{???}
N(k)\le \tfrac{\bigl(k! (k-1)!\bigr)^2}{(2\pi)^2 x_1^k x_2^{k-1}}\int\limits_{\xi_s\in (-\pi,\pi)}
\prod_{j\ge 0}\Big|H\Bigl(\tfrac{z_1z_2^j}{(j!)^2}\Bigr)\Big|\,d\xi_1 d\xi_2,\,\,z_s=x_s e^{i\xi_s}.
\end{equation}

\begin{lemma}\label{lem1} %For $z=|z|e^{i\theta}$, $\theta\in (-\pi,\pi)$, 
For $H(z)=\sum_{r\ge 0}\tfrac{z^r}{(r!)^2}$,  $\text{arg}(z)\in (-\pi,\pi]$, we have 
\begin{equation*}%\label{3.12}
|H(z)|\le \tfrac{|\exp(2z^{1/2})|}{\max(1,\a |z|^{1/4})},\quad z^{1/2}:=|z|^{1/2} \exp\bigl(\tfrac{1}{2}\text{arg}(z)\bigr),
%H(|z|) \exp\bigl[-|z|(1-\cos(\text{arg}(z))\bigr]=H(|z|) \exp\bigl[-(|z|-\text{Re}\,(z))\bigr].
\end{equation*}
%holds.
where $\a>0$ is an absolute constant.
\end{lemma}
\si
%Visually and technically,  the inequality to prove is similar to these two inequalities
%\begin{align*}
%&\Big|\sum_{j\ge a}\tfrac{z^j}{j!}\Bigr|\le \Bigl(\sum_{j\ge a}\tfrac{|z|^j}{j!}\Bigr)\cdot\exp\bigl(-\tfrac{|z|(1-\cos\theta)}{a+1}\bigr),\\
%&\Big|\tfrac{1}{1-z}\Big|\le \tfrac{1}{1-|z|}\cdot\exp\bigl(-|z|(1-\cos\theta)\bigr),\quad |z|<1,
%\end{align*}
%see \cite{Pit1} and \cite{Pit2} respectively.
%\si
\begin{proof} Introduce $I(z):=\sum_{r=0}^{\infty}\tfrac{z^{2r}}{(r!)^2}$, $\text{arg}(z)\in(-\pi/2,\pi/2]$, so that $H(z^2)=I(z)$.  It suffices to prove that, for $\arg(z)\in (-\pi/2,\pi/2]$,
\begin{equation}\label{2}
|I_0(2z)|\le \tfrac{\big|\exp(2z)\big|}{\max(1,\a |z|^{1/2})},\quad I_0(\xi)=\sum_{j\ge 0}\frac{\bigl(\tfrac{\xi^2}{4}\bigr)^j}{(j!)^2},
\end{equation}
where $I_0(\xi)$ is the zero-order, modified, Bessel function, since $H(z^2)=I_0(2z)$. Here is a surprisingly simple integral formula for $I_0(2z)$ that does the job:
\begin{equation}\label{2.1}
I_0(2z)=\pi^{-1}\exp(2z)\int_0^4\tfrac{\exp(-zy)}{\sqrt{y(4-y)}}\,dy,\quad \text{arg}(z)\in (-\pi/2,\pi/2].
\end{equation}
The equation \eqref{2.1} follows from  identity 10.32.2 in \cite{DLMF}, that is, 
\[
I_0(2z)=\tfrac{1}{\pi}\int_1^1(1-t^2)^{-1/2}e^{-2tz}\,dt,
\]
via substitution $y=2t+2$. (Our original proof of \eqref{2.1} was from scratch. We are grateful to the reviewer for the above reference.) Define
    \begin{equation}
    \label{eq:lap1}
 G(z)=\int\limits_0^4\!\!\frac{e^{-zy}}{\sqrt{y(4-y)}}dy,  \end{equation}
so that $I_0(2z)=\pi^{-1}\exp(2z)G(z)$, whence $|I_0(2z)|=\pi^{-1}|\exp(2z)| \cdot |G(z)|$. Here for $\Re z\ge 0$
\begin{equation*}
|G(z)|\le\int_0^4\tfrac{1}{\sqrt{y(4-y)}}\,dy=\pi.
\end{equation*}
So, by \eqref{2.1},  $|I_0(2z)|\le \pi^{-1}|\exp(2z)|\pi=|\exp(2z)|$. Next, 
let us show that for $g(z):=|G(z)|/G(|z|)$ we have $\sup_{\mathbb{H}}g(z)<\infty$,
where $\mathbb{H}:=\{z: \Re z\ge 0\}$.
The function $g$ is well-defined and continuous in $\mathbb{H}$, since $G(z)$ is continuous in $\mathbb{H}$ and $G(|z|)>0$.    %$\displaystyle \sup_{z\in\mathbb{H}} g(z)\le \max\{2^{9/4},\max_{|z|\le 1} g(z)\}$\footnote{A concrete bound for $g$ for $|z|\le 1$ can be easily obtained using the Maclaurin expansion of $I_0(2z)$.}
%\end{Lemma}
%\begin{proof}
 %Continuity follows from the continuity of $G$ and positivity of $G(|z|)$ in $\mathbb{H}$. 
We are left with showing $\sup_{z\in \mathbb{H}: |z|\ge 1}g(z)<\infty$.
%for $|z|\ge 1$, $z\in\mathbb{H}$ which we prove next.
Noting that $G(\overline{z})=\overline{G(z) }$, it suffices to show this boundedness in the closed fourth quadrant $Q_{IV}$. For such $z$ we have
  \begin{equation}
    \label{eq:g-homot}
    G(z)=\int_{\mathcal{C}_a} \frac{e^{-zy}}{\sqrt{y(4-y)}}dy,
  \end{equation}
  where $a>0$ and $\mathcal{C}_a$ is the oriented polygonal line joining the points $y_1=0,  y_2=a\sqrt{i},  y_3=a\sqrt{i}+4$, and $y_4=4$. %Observe that the (absent) line segment $(y_2,y_4)$ is perpendicular to the segment $(y_1,y_2)$. 
Since $|e^{-zy}|\le 1$ in $\mathbb{H}$ and $\min\{|y|,|4-y|\}\ge \Im y= 2^{-1/2} a$, the integral over the segment $(y_2,y_3)$ is bounded by 
$4\sqrt{2} /a$,
 hence it goes to zero as $a\to \infty$. Passing to the limit $a\to \infty$, we thus get
  \begin{equation}
    \label{eq:g-homot2}
    \begin{aligned}
    G(z)&=\int_{0}^{\infty e^{i\pi/4}} \frac{e^{-zy}}{\sqrt{y(4-y)}}dy-\int_{4}^{\infty e^{i\pi/4}+4} \frac{e^{-zy}}{\sqrt{y(4-y)}}dy\\
    &=\int_{0}^{\infty e^{i\pi/4}}\left(\frac{e^{-zu}}{\sqrt{u(4-u)}}-i\frac{e^{-4z-zu}}{\sqrt{(4+u)u}}\right)du.
  \end{aligned}
  \end{equation}
  In the last integral  we have $|\sqrt{u(4-u)}|\ge 2^{3/4}\sqrt{|u|}$, $|\sqrt{u(4+u)}|\ge 2 \sqrt{|u|}$, $|e^{-4z}|\le 1$ and (since $z\in Q_{IV}$ and $\arg u=\frac{\pi}{4}$), $|e^{-zu}|\le e^{-2^{-1/2}|z||u|}$.  Since $2^{-3/4}+2^{-1}<2$, we get
  \begin{equation}
    \label{eq:ineq1}
    |G(z)|\le 2\int_0^\infty u^{-1/2} e^{-2^{-1/2}|z|u}\, d u=2^{5/4}\sqrt{\pi} |z|^{-1/2}.
  \end{equation}
Importantly, we need this bound only for $|z|\ge 1$. For those $z$'s, we have
%On the other hand, for  $|z|\ge 1$%, $\exp(-4|z|)<1/2$ and
  \begin{equation}
    \label{eq:gabs}
    G(|z|)\ge\tfrac{1}{2} \int_0^4\frac{e^{-|z|y}}{\sqrt{y}}dy\ge \tfrac{1}{2}|z|^{-1/2}\int_0^4\eta^{-1/2}e^{-\eta}\,d\eta.\\
  \end{equation}  
    %=\tfrac{1}{2}\int_0^\infty\frac{e^{-|z|y}}{\sqrt{y}}dy-\tfrac{1}{2}\int_4^{\infty}\frac{e^{-|z|y}}{\sqrt{y}}dy\\
   % \ge \frac12 \int_0^\infty\frac{e^{-|z|y}}{\sqrt{y}}dy=\frac{\sqrt{\pi}}{2\sqrt{|z|}}.
  %\end{multline}
Combining \eqref{eq:ineq1} and \eqref{eq:gabs}, we see that $\sup_{z\in \Bbb H:|z|\ge 1}\tfrac{|G(z)|}{G(|z|)}<\infty$.
Therefore
\[
|I_0(2z)|=O\bigl(|\exp(2z)| \cdot G(|z|)\bigr)=O\bigl(|\exp(2z)|\cdot |z|^{-1/2}\bigr),
\]
uniformly for $z\in\{ \Bbb H:|z|\ge 1\}$. This bound and the inequality $|I_0(2z)|\le |\exp(2z)|$, $\Re z\ge 0$, complete 
the proof of the lemma.
% \end{proof} 
\end{proof}
%{\bf Note.\/}  A closer look at $G(z)$ (see Appendix) shows that for $\text{Re}(z)\ge 0$, $G(|z|)=\Theta\bigl((1+|z|)^{-1/2}
%\bigr)z|)$ and $G(z)/G(|z|)$ is bounded. So,
%\[
%|I(2z)|=O\Bigl(\tfrac{|\exp(2z)|}{1+|z|^{1/2}}\Bigr),\quad \text{Re}(z)\ge 0,
%\]
%which is qualitatively the best upper bound for $|I(2z)|$.
\subsection{Proof of Theorem \ref{thmStar}, part {\bf(a)}} We use the lemma \ref{lem1} to prove that, for 
\[
\delta>0,\quad K_n=\big\lfloor\exp\bigl((2+\delta)\sqrt{\log n \cdot\log(\log n)}\bigr)\big\rfloor, 
\]
we have $\sum_{k\ge K_n}m_n(k)\to 0$, where $m_n(k)$ is the expected number of pairs of twins of size $k$.

We turn back to the integrand in \eqref{???}. With $z_s=x_s e^{i\xi_s}$,
and $X_j:=\tfrac{x_1x_2^j}{(j!)^2}$, by Lemma \ref{lem1} and \eqref{3.21}, we have
\begin{equation}\label{3.02}
\prod_{j\ge 0}\Big|H\Bigl(\tfrac{z_1z_2^j}{(j!)^2}\Bigr)\Big|\le \prod_{j\ge 0}\tfrac{\exp(2X_j^{1/2})}{\max(1,\a |X_j|^{1/4})}
\exp\Bigl(X_j\bigl(\cos(\xi_1+j\xi_2)-1\bigr))\Bigr).
\end{equation}
%H\Bigl(\tfrac{x_1x_2^j}{(j!)^2}\Bigr)
%\exp\biggl(\tfrac{x_1x_2^j}{(j!)^2}\bigl(\cos(\xi_1+j\xi_2)-1\bigr)\biggr).
%\end{equation}
%\le \exp\Bigl(2\tfrac{x_1^{1/2}x_2^{j/2}}{j!}\Bigr)
%\exp\biggl(\tfrac{x_1x_2^j}{(j!)^2}\bigl(\cos(\xi_1+j\xi_2)-1\bigr)\biggr).
%\end{multline*}

Let us look closely at $\{X_j\}$. Using $j!\le\bigl(\tfrac{j+1}{2}\bigr)^j$, $j\ge 1$, we have
\[
X_j\ge x_1\cdot\tfrac{x_2^j}{\bigl(\tfrac{j+1}{2}\bigr)^{2j}}=x_1\exp\Bigl(2 j\log\tfrac{2x_2^{1/2}}{j+1}\Bigr)
=x_1\exp\Bigl(2 j\log\tfrac{2\tfrac{k-1}{k}}{j+1}\Bigr).
\]
The function  $\eta \log\tfrac{2\tfrac{k-1}{k}}{\eta+1}$ is decreasing for $\eta\ge 1$. Pick $\eps\in (0,1/2)$ and introduce
$J=\Big\lfloor\eps\tfrac{\log x_1}{\log(\log x_1)}\Big\rfloor$. Then, uniformly for $j\le J$, 
\begin{multline}\label{3.031}
X_j\ge x_1\exp\Bigl(2 J\log\tfrac{2\tfrac{k-1}{k}}{J+1}\Bigr)\\
=x_1\exp(-2J\log J +O(J))=\exp\bigl((1-2\eps)\log x_1+o(\log x_1)\bigr)\to\infty.
\end{multline}
Further
\[
\sum_{j=1}^{J-1}2 j\log\tfrac{2\tfrac{k-1}{k}}{j+1}\ge 2\int_1^{J} \eta\log\tfrac{2\tfrac{k-1}{k}}{\eta+1}\,d\eta\ge -J^2\log J;
\]
the last inequality follows easily from integration by parts. 
Therefore
\begin{equation}\label{3.0315}
\begin{aligned}
\prod_{j<J}X_j&\ge x_1^J\exp(-J^2\log J)\ge\exp\bigl[J(\log x_1-J\log J)\bigr]\\
&\ge \exp\Bigl(\eps'(1-\eps')\tfrac{\log^2 x_1}{\log(\log x_1)}\Bigr)\\
&=x_1^{\eps'(1-\eps')\tfrac{\log x_1}{\log(\log x_1)}},\quad\forall \eps'<\eps.
\end{aligned}
\end{equation}
The point of this bound is that the product in question exceeds $x_1$ raised to a sub-logarithmic power, still approaching infinity when $k$ does. 
%Furthermore, by \eqref{3.031},
%\begin{equation*}
%\sum_{j<J}X_j^{-1/2}\le J\exp\bigl(-\tfrac{(1-2\eps)\log x_1}{2}+o(\log x_1)\bigr)\le 1,
%\end{equation*}
%for large $k$, since $\log J=o(\log x_1)$. 
It follows from \eqref{3.0315} that
\begin{multline*}
 \prod_{j\le J}\tfrac{\exp(2X_j^{1/2})}{\max(1,\a |X_j|^{1/4})}
\le  \exp\Bigl(-\eps'(1-\eps')\tfrac{\log^2 x_1}{4\log(\log x_1)}+J\log\a\Bigr)\exp\Bigl(2\sum_{j\le J}X_j^{1/2}\Bigr)\\
= \exp\Bigl(-\eps'(1-\eps')\tfrac{\log^2 x_1}{4\log(\log x_1)}+O(\eps\log x_1)\Bigr)\exp\Bigl(2\sum_{j<J}X_j^{1/2}\Bigr).
\end{multline*}
Since $\tfrac{\exp(2X_j^{1/2})}{\max(1,\a |X_j|^{1/4})}\le \exp(2X_j^{1/2})$ for all $j$, the above inequality implies that for $\eps''<\eps'$,
\begin{multline}\label{3.032}
\prod_{j\ge 0}\tfrac{\exp(2X_j^{1/2})}{\max(1,\a |X_j|^{1/4})}
\le  \exp\Bigl(-\tfrac{\eps''(1-\eps'')\log^2 x_1}{4\log(\log x_1)}\Bigr)\exp\Bigl(2\sum_{j\ge 0}X_j^{1/2}\Bigr)\\
= \exp\Bigl(-\tfrac{\eps''(1-\eps'')\log^2 x_1}{4\log(\log x_1)}\Bigr)\exp\bigl(2x_1^{1/2}\exp(x_2^{1/2})\bigr)\\
\le e^{2k}\cdot \exp\Bigl(-\tfrac{\eps{''}(1-\eps{''})\log^2k}{\log(\log k)}\Bigr),%,\quad\forall\,\eps^{''}<\eps'<\eps,
\end{multline}
as  $x_1=k^2\exp\bigl(-\tfrac{2(k-1)}{k}\bigr)$, $x_2=\bigl(\tfrac{k-1}{k}\bigr)^2$.

Dropping the factors $\exp\Bigl(X_j\bigl(\cos(\xi_1+j\xi_2)-1\bigr))\Bigr)(\le 1)$ with $j>1$ in \eqref{3.02}, and using
\eqref{3.032}, we obtain
\begin{multline}\label{3.115}
\prod_{j\ge 0}\Big|H\Bigl(\tfrac{z_1z_2^j}{(j!)^2}\Bigr)\Big|\le e^{2k}\cdot \exp\Bigl(-\tfrac{\eps{''}(1-\eps^{''})\log^2k}{\log(\log k)}\Bigr)
\\
\times\exp\Bigl(x_1 (\cos\xi_1-1)+x_1x_2\bigl(\cos(\xi_1+\xi_2)-1\bigr)\Bigr),\quad \xi_j\in (-\pi,\pi).
\end{multline}
In order to use this last estimate on the right-hand side of \eqref{???}, we need to find an upper bound for the {\em integral} of the exponential factor above.
Let $I=(-\pi,\pi)$. In order to obtain the needed upper bound,  let us set
\[f(\xi_1,\xi_2):=x_1(\cos (\xi_1)-1) +x_1x_2(\cos(\xi_1+\xi_2)-1).\]
Clearly, $f(\xi_1,\xi_2)\le 0$, with equality holding if and only if 
$\cos(\xi_1)=1$,  and  $\cos(\xi_1+\xi_2)=1$, where $\xi_i\in I$,
or equivalently if and only if  $\xi_1=\xi_2=0$. 

For $\max \{|\xi_1|,|\xi_2|\} \geq \pi/2$, we have
\begin{equation} \label{star} f(\xi_1,\xi_2)\le -ak^2< 0,  \end{equation}
for some positive constant $a$. For  $\max \{|\xi_1|,|\xi_2|\} \leq \pi/2$,
using the Taylor expansion of $\cos$, we obtain
\begin{equation} \label{twostar} f(\xi_1,\xi_2)\le -b [ x_1 \xi_1^2 + x_1x_2(\xi_1+\xi_2)^2]\le -c k^2 [\xi_1^2 +\xi_2^2] .\end{equation}
Combining \eqref{star} and \eqref{twostar},  we get
\begin{equation}\label{3star}
f(\xi_1,\xi_2)\le -d k^2 [\xi_1^2 +\xi_2^2], 
\end{equation}
for $\xi_i\in I$, for $ i=1,2$, and some fixed constant $d>0$. 

Therefore, %Now we are in a position to prove an upper bound for $\int_{I\times I} \exp(f(\xi_1,\xi_2)) \ d\xi_1 d\xi_2$ as follows.
\begin{multline}\label{4star}
\int_{I\times I}\! \exp(f(\xi_1,\xi_2)) \ d\xi_1 d\xi_2\\
\le\int\limits_{I\times I}\!\! \exp[ x_1(\cos(\xi_1-1)+x_1x_2(\cos(\xi_1+\xi_2)-1)] \,d\xi_1d \xi_2\\
\le\int\limits_{I\times I} \!\!\exp(-dk^2(\xi_1^2+\xi_2^2)) \ d\xi_1 d\xi_2\\ 
 =\! \int\limits_{I} \!\exp(-dk^2 \xi_1^2) \,d\xi_1 \cdot   \int\limits_{I}\! \exp(-dk^2 \xi_2^2) \,d\xi_2
\leq \frac{\pi}{dk^2}.%=O(1) .
\end{multline}
%The bound is crude, but it is all we need, since the second exponential factor on the top RHS of \eqref{3.115} is much smaller. And it is this factor that makes what comes next possible. 
Putting together \eqref{???}, \eqref{3.115} and the last inequality, %and using $x_1=k^2\exp\bigl(-\tfrac{2(k-1)}{k}\bigr)$,
%$x_2=\bigl(\tfrac{k-1}{k}\bigr)^2$, 
we conclude that  
\begin{multline*}
N(k)=O\Bigl(\frac{(k!)^4e^{2k}}{x_1^k x_2^{k-1}k^4}\Bigr)\exp\Bigl(-\tfrac{\eps{''}(1-\eps{''})\log^2k}{\log(\log k)}\Bigr)\\
=O\Bigl(\tfrac{(k!)^4e^{4k}}{k^{2k+4}}\Bigr)\exp\Bigl(-\tfrac{\eps{''}(1-\eps{''})\log^2k}{\log(\log k)}\Bigr).
\end{multline*}
This bound and \eqref{4} imply that
\begin{multline*}
S_n(k) = O\left(\binom{n}{k}\binom{n-k}{k}(n-2k)^{n-2k+1}N(k)\right)\\
=O\Bigl(\tfrac{n!(n-2k)^{n-2k+1}(k!)^2 k^{-2k-4}e^{4k}}{(n-2k)!}\Bigr)\exp\Bigl(-\tfrac{\eps{''}(1-\eps{''})\log^2k}{\log(\log k)}\Bigr).\\
%=O\left((k!)^4k^{-2k-2}e^{4k}\right) \exp\Bigl(-\tfrac{\eps{''}(1-\eps^{''})\log^2k}{\log(\log k)}\Bigr),\quad\forall \eps^{''}<\eps.
\end{multline*}
Using Stirling's formula for factorials we obtain that 
\begin{align*}
S_n(k)&= O\left(\binom{n}{k}\binom{n-k}{k}(n-2k)^{n-2k+1}N(k)\right)\exp\Bigl(-\tfrac{\eps{''}(1-\eps{''})\log^2k}{\log(\log k)}\Bigr)\\
&=O\Bigl(\tfrac{n^{n+1}}{k^3}\Bigr)\exp\Bigl(-\tfrac{\eps{''}(1-\eps{''})\log^2k}{\log(\log k)}\Bigr).
\end{align*}
%O\biggl(n^{n+1} \exp\Bigl(-\tfrac{\eps{''}(1-\eps{''})\log^2k}{\log(\log k)}\Bigr)\biggr).
%\end{align*}
%& O\biggl( \frac{n!}{k!^2(n-2k)!} \cdot (n-2k)^{n-2k+1} \cdot   \frac{\bigl(k!(k-1)!\bigr)^2e^{2k}}{k^{2k}e^{-2k}k^{5/2}}\biggr) \\
%& = & O\biggl(\frac{n^{n+1}}{  e^{2k} }\cdot   \frac{(k-1)^{2k-1}}{k^{2k}} \cdot \frac{e^{2k}}{e^{-2k} k^{5/2}}     \biggr) \\
%& = & O\biggl(\frac{n^{n+1}}{k^{7/2}}\biggr).
%\end{eqnarray*}
Dividing this bound by the number $n^{n-1}$ of all rooted Cayley trees, we get that the expected number of the pairs of twins of size $k$, i.e. $m_n(k)=\tfrac{S_n(k)}{n^{n-1}}$, is
\begin{equation}\label{3.05}
m_n(k)=O\biggl(\tfrac{n^2}{k^3}\exp\Bigl(-\tfrac{\eps{''}(1-\eps{''})\log^2k}{\log(\log k)}\Bigr)\biggr),\quad\forall\eps''<\eps<1/2.
\end{equation}
Since the series $\sum_{k\ge 1} k^{-3}$ converges, it follows that for \newline  $k_n=\exp\Bigl((2+\delta)\sqrt{\log n \cdot\log(\log n)}\biggr)$ and $\delta>0$ we have
\begin{multline}\label{3.06}
\sum_{k\ge k_n}m_n(k)=O\Bigl(n^2\exp\Bigl(-\tfrac{\eps{''}(1-\eps{''})\log^2k_n}{\log(\log k_n)}\Bigr)\Bigr)\sum_{k\ge 1}k^{-3}\\
=O\exp\Bigl(2\bigl[(1-(2+\delta)^2\eps''(1-\eps''))\bigr]\log n+o(\log n)\Bigr),
\end{multline}
uniformly for all $\eps''<1/2$. Picking $\eps''$ sufficiently close to $1/2$, we conclude that 
$\lim_{n\to\infty}\sum_{k\ge k_n}m_n(k)=0$.

\subsection{Proof of Theorem \ref{thmStar}, part {\bf (b)\/}} The claim is: If $\delta\in\! (0,2)$, $k_n\!:=\!\Big\lfloor\!\exp\Bigl(\!(2-\delta)\sqrt{\log n\cdot \log(\log n)}\Bigr)\!\Big\rfloor$, 
then $m_n(k_n)\to\infty$.

%{\bf Note.\/} We are tempted to conjecture that, with probability approaching $1$, the largest size of twin subtrees is  of order $\exp\Bigl(\Theta\bigl(\!\sqrt{\log n\cdot \log(\log n)}\,\bigr)\Bigr)$. Could it be that, in fact, the distribution of the largest twin
%size is concentrated around $\big\lfloor\exp\Bigl(2\sqrt{2}\sqrt{\log n\cdot \log(\log n)}\Bigr)\big\rfloor$?
\begin{proof} From here on we can, and will consider the twins of maximum out-degree $<d=d_n$, with $d_n$ yet 
to be chosen. Let us upper-bound the contribution to the RHS of \eqref{???} for $k=k_n$, coming from $(z_1=x_1e^{i\xi_1}, z_2=x_2e^{i\xi_2})$, where 
$x_1=k^2\exp\bigl(-\tfrac{2(k-1)}{k}\bigr)$, $x_2=\bigl(\tfrac{k-1}{k}\bigr)^2$,
and also $\|\boldsymbol\xi\|\le k^{-1/2+\eps}$, $\eps\in (0,1/2)$. By Lemma \ref{lem1},
\begin{multline*}
\prod_{j<d}\Big|H\Bigl(\tfrac{z_1z_2^j}{(j!)^2}\Bigr)\Big|\le\exp\biggl(2x_1^{1/2}\sum_{j\le d}\tfrac{x_2^{j/2}}{j!}\Re\bigl(\exp(i(\xi_1/2+j\xi_2/2)\bigr)\biggr)\\
\le\exp\biggl(2x_1^{1/2}\sum_{j<d}\tfrac{x_2^{j/2}}{j!}\Re\bigl(\exp(i(\xi_1/2+j\xi_2/2)\bigr)+O(x_1^{1/2}/d!)\biggr)\\
=\exp\bigl(2x_1^{1/2}W(\xi_1,\xi_2)+O(x_1^{1/2}x_2^{d/2}/d!)\bigr),
\end{multline*}
where
\[
W(\xi_1,\xi_2):=\cos\bigl(\tfrac{\xi_1}{2}+x_2^{1/2}\sin\tfrac{\xi_2}{2}\bigr)\exp\bigl(x_2^{1/2}\cos\tfrac{\xi_2}{2}\bigr).
\]
%Recall that $x_1=k^2\exp\bigl(-\tfrac{2(k-1)}{k}\bigr)\sim k^2e^{-2}$, $x_2=\bigl(\tfrac{k-1}{k}\bigr)^2\sim 1$. In particular, 
Since $x_1=\Theta(k^2)$, $x_2=\Theta(1)$,
the big-O term is $o(1)$, if $d= \lfloor2\log k/\log\log k\rfloor$, which we assume from now.

Clearly, $W(\xi_1,\xi_2)\le \exp(x_2^{1/2})=W(0,0)$, and it is easy to check that $(0,0)$ is a single stationary point of
$W(\xi_1,\xi_2)$ in $[-\pi,\pi]^2$, whence it is a unique maximum point of $W(\xi_1,\xi_2)$ in this square. In addition,
\begin{equation*}
W''_{\xi_1}(0,0)=-\tfrac{e^{x_2^{1/2}}}{4},\,\,W''_{\xi_2}(0,0)=-\tfrac{(x_2+x_2^{1/2})e^{x_2^{1/2}}}{4},\,\,
W''_{\xi_1,\xi_2}(0,0)=-\tfrac{x_2^{1/2}e^{x_2^{1/2}}}{4},
\end{equation*}
so that $W''_{\xi_1}(0,0), W''_{\xi_2}(0,0)<0$, and
\[
W''_{\xi_1}(0,0)\cdot W''_{\xi_2}(0,0)-\bigl(W''_{\xi_1,\xi_2}(0,0)\bigr)^2=\tfrac{x_2^{1/2}e^{2x_2^{1/2}}}{16}>0.
\]
This implies that $W(\xi_1,\xi_2)$ is strictly concave in the vicinity of $(0,0)$, and moreover 
$W(\xi_1,\xi_2)\le W(0,0)-\be(\xi_1^2+\xi_2^2)$, $\xi_s\in [-\pi,\pi]$, for a constant $\be>0$ as $k\to\infty$ because 
$x_2=\Theta(1)$. So, 
\[
\prod_{j<d}\Big|H\Bigl(\tfrac{z_1z_2^j}{(j!)^2}\Bigr)\Big|\le \exp\bigl(2x_1^{1/2}e^{x_2^{1/2}}
-\Theta(k\|\boldsymbol\xi\|^2)+o(1)\bigr).
\]
Therefore 
\begin{multline}\label{onestar}
\int\limits_{\|\bold \xi\|\ge k^{-1/2+\eps}}\tfrac{1}{z_1^{k+1} z_2^k}\prod_{j< d}H\Bigl(\tfrac{z_1z_2^j}{(j!)^2}\Bigr)\,dz_1 dz_2\\
 =O\biggl(\!\tfrac{\exp\bigl(2x_1^{1/2}e^{x_2^{1/2}}\bigr)}{x_1^kx_2^{k-1}}\!\!\!\!\!\!\!\int\limits_{\|\boldsymbol\xi\|\ge k^{-1/2+\eps}}\!\!\!\!\!\!\!\!\!     e^{-\Theta(k\|\bold\xi\|^2)}\,d\xi_1 d\xi_2\!\biggr)=O\biggl(\!\tfrac{\exp\bigl(2k-\Theta(k^{2\eps})\bigr)}{x_1^kx_2^{k-1}}\biggr).
\end{multline}
It remains to sharply evaluate the contribution to the Cauchy integral coming from $\boldsymbol\xi$'s with $\|\boldsymbol\xi\|\le k^{-1/2+\eps}$. By Lemma \ref{lem1}, we have 
\[
H(z)=\pi^{-1}\exp(2\sqrt{z})\int_0^4\tfrac{\exp(-\sqrt{z} y)}{\sqrt{y(4-y)}}\,dy,\quad \text{arg}(z)\in [-\pi,\pi].
\]
And, for the values of  $z_1, z_2$, and  $j$ in question, $\text{arg}(z_1z_2^j)=\xi_1+j\xi_2=O(k^{-1/2+\eps}d)=o(1)$. Consequently,
\begin{equation}\label{3.055}
\tfrac{1}{z_1^k z_2^{k-1}}\!\prod_{j< d}\!H\Bigl(\!\tfrac{z_1z_2^j}{(j!)^2}\!\Bigr)\!=\!\tfrac{\pi^{-d}}{z_1^k z_2^{k-1}}\prod_{j<d}\!\exp\bigl(\tfrac{2z_1^{1/2}z_2^{j/2}}{j!}\bigr)
\!\!\int_0^4\!\tfrac{\exp\bigl(-\tfrac{z_1^{1/2}z_2^{j/2}}{j!} y_j\bigr)}{\sqrt{y_j(4-y_j)}}\!dy_j.
\end{equation}
Here
\begin{equation}\label{3.06}
\tfrac{1}{z_1^k z_2^{k-1}}\prod_{j<d}\exp\bigl(\tfrac{2z_1^{1/2}z_2^{j/2}}{j!}\bigr)=
\tfrac{\exp\bigl(2z_1^{1/2}e^{z_2^{1/2}}\bigr)+O(x_2^{d/2}/d!)}{z_1^kz_2^{k-1}},
\end{equation}
and 
\begin{multline*}
2z_1^{1/2}e^{z_2^{1/2}}=2x_1^{1/2}\exp(i\xi_1/2+x_2^{1/2}e^{i\xi_2/2})\\
=2x_1^{1/2}e^{x_2^{1/2}}\exp\bigl(i\xi_1/2+x_2^{1/2}(e^{i\xi_2/2}-1)\bigr)\\
=2x_1^{1/2}e^{x_2^{1/2}}\bigl[1+(i\xi_1/2+x_2^{1/2}(e^{i\xi_2/2}-1))+\tfrac{1}{2}(i\xi_1/2+x_2^{1/2}(e^{i\xi_2/2}-1))^2+
O(\|\boldsymbol\xi\|^3)\bigr]\\
=2x_1^{1/2}e^{x_2^{1/2}}\bigl[1+i(\xi_1/2+x_2^{1/2}\xi_2/2)\\
 -\xi_1^2/8-(x_2+x_2^{1/2})\xi_2^2/8-\xi_1\xi_2 x_2^{1/2}/4+O(\|\boldsymbol\xi\|^3)\bigr].
\end{multline*}
Therefore
\begin{multline*}
2z_1^{1/2}e^{z_2^{1/2}}-k\log z_1-(k-1)\log z_2=2x_1^{1/2}e^{x_2^{1/2}}-k\log x_1-(k-1)\log x_2\\
+i\xi_1\bigl(x_1^{1/2}e^{x_2^{1/2}}-k\bigr)+i\xi_2(x_1^{1/2}e^{x_2^{1/2}}x_2^{1/2}-(k-1)\bigr)\\
-\tfrac{1}{4}x_1^{1/2}e^{x_2^{1/2}}\bigl(\xi_1^2+(x_2+x_2^{1/2})\xi_2^2-2\xi_1\xi_2x_2^{1/2}\bigr)+O(k\|\boldsymbol\xi\|^3).
\end{multline*}
Recalling that $\|\boldsymbol\xi\|\le k^{-1/\eps}$, we see that $O(k\|\boldsymbol\xi\|^3)=O(k^{-1/2+3\eps})\to 0$, if $\eps<1/6$, which we assume from now. Furthermore, the linear combination of $i\xi_1$ and $i\xi_2$ in the above sum disappears, thanks to
\[
x_1^{1/2}e^{x_2^{1/2}}-k=0,\quad x_1^{1/2}e^{x_2^{1/2}}x_2^{1/2}-(k-1)=0,
\]
the conditions we were led to in our preliminary attempt to bound $N(k)$.
%meaning that 
%\[
%x_1=\bigl(\tfrac{k}{2}\bigr)^2\exp\bigl(-\bigl(\tfrac{k-1}{k}\bigr)^2\bigr),\quad x_2=\bigl(\tfrac{k-1}{k}\bigr)^2.
%\]
%meeting our constraint $x_1=\Theta(k^2)$, $x_2=O(1)$. For this choice of $x_1$ and $x_2$, 
So, we have
\begin{multline*}
2z_1^{1/2}e^{z_2^{1/2}}-k\log z_1-(k-1)\log z_2=2k-k\log x_1-(k-1)\log x_2
\\
-\tfrac{k}{4}\bigl(\xi_1^2+(x_2+x_2^{1/2})\xi_2^2-2\xi_1\xi_2x_2^{1/2}\bigr)+O(k\|\boldsymbol\xi\|^3).
\end{multline*}
Therefore,  equation \eqref{3.06} becomes 
\begin{multline}\label{3.07}
\tfrac{1}{z_1^k z_2^{k-1}}\prod_{j<d}\exp\bigl(\tfrac{2z_1^{1/2}z_2^{j/2}}{j!}\bigr)\\
=\tfrac{e^{2k}}{x_1^k x_2^{k-1}}\exp\bigl[-\tfrac{k}{4}\bigl(\xi_1^2+(x_2+x_2^{1/2})\xi_2^2-2\xi_1\xi_2x_2^{1/2}\bigr)+O(k\|\boldsymbol\xi\|^3+e^{-k})\bigr].\\
%\cdot \bigl(1+O(1/(d! e^{k}))\bigl).
\end{multline}
The quadratic form is negative definite, since $x_2=\Theta(1)$. We still have to evaluate 
\begin{equation}\label{3.08}
\begin{aligned}
&\qquad\qquad\qquad\quad\prod_{j<d}\int_0^4\tfrac{\exp\bigl(-Z_j y_j\bigr)}{\sqrt{y_j(4-y_j)}}\,dy_j,\\
&Z_j=X_j \exp(i(\xi_1/2+j\xi_2/2)),\quad
X_j:=\tfrac{k\exp\bigl(-\tfrac{1}{2}\bigl(\tfrac{k-1}{k}\bigr)^2\bigr)}{j!}\bigl(\tfrac{k-1}{k}\bigr)^j,\\
\end{aligned}
\end{equation}
where $d= \lfloor2\log k/\log\log k\rfloor$, $\|\boldsymbol\xi\|\le k^{-1/2+\eps}$, $\eps<1/6$. 

{\bf (a)\/} Suppose that $j< j(k):=\lfloor\tfrac{\log k}{\log\log k}\rfloor$. Then 
\[
\tfrac{k}{j!}\ge \frac{k}{j(k)!}=\exp\Bigl(\Theta\bigl(\tfrac{\log k\,\cdot\,\log\log\log k}{\log\log k}\bigr)\Bigr)\to\infty.
\]
So, using
\[
I_0(z)\sim\tfrac{e^z}{(2\pi z)^{1/2}}\bigl(1+O(z^{-1})\bigr),\quad |\text{arg}(z)|\le \tfrac{\pi}{2}-\delta,\,\,\delta\in (0,\pi/2),
\]
(see \cite{DLMF}), we have	
\begin{multline*}
\int_0^4\tfrac{\exp\bigl(-Z_j y_j\bigr)}{\sqrt{y_j(4-y_j)}}\,dy_j=\pi \exp(-2Z_j)I_0(2Z_j)\\
=\pi\exp(-2Z_j)\tfrac{e^{2Z_j}}{(4\pi Z_j)^{1/2}}\bigl(1+O(Z_j^{-1})\bigr)=\tfrac{\pi}{(4\pi Z_j)^{1/2}}\bigl(1+O(Z_j^{-1})\bigr)\\
=\tfrac{\pi}{(4\pi X_j)^{1/2}}\cdot\bigl[1+O\bigl(\exp\bigl(-\Theta\bigl(\tfrac{\log k\,\cdot\,\log\log\log k}{\log\log k}\bigr)\bigr)\bigr)\bigr],
%=\int_0^4\tfrac{\exp\bigl(-X_j y_j\bigr)}{\sqrt{y_j(4-y_j)}}\,dy_j\cdot\bigl[1+O\bigl(\exp\bigl(-\Theta\bigl(\tfrac{\log k\,\cdot\,\log\log\log k}{\log\log k}\bigr)\bigr)\bigr)\bigr],
\end{multline*}
implying that
\begin{align*}
\prod_{j< j(k)}\int_0^4\tfrac{\exp\bigl(-Z_j y_j\bigr)}{\sqrt{y_j(4-y_j)}}\,dy_j
&=(1+O(\Delta_k))\prod_{j<j(k)}\tfrac{\pi}{(4\pi X_j)^{1/2}},\\
%\int_0^4\tfrac{\exp\bigl(-X_j y_j\bigr)}{\sqrt{y_j(4-y_j)}}\,dy_j ,\\
\Delta_k&:=\tfrac{\log k}{\log\log k}\exp\bigl(-\Theta\bigl(\tfrac{\log k\,\cdot\,\log\log\log k}{\log\log k}\bigr)\bigr).
\end{align*}
Here, by \eqref{3.08}, $X_j\le \tfrac{k}{j!}$. So, using Stirling's formula with the remainder term, namely
$
j!=(2\pi j)^{1/2}\bigl(\tfrac{j}{e}\bigr)^j [1+O(1/(j+1))],
$
we obtain
\begin{multline}\label{3.09}
\prod_{j<j(k)}\tfrac{\pi}{(4\pi X_j)^{1/2}}=\exp[O(j(k)\log j(k))]\cdot k^{-j(k)/2}\prod_{j<j(k)}\bigl(\tfrac{j}{e}\bigr)^{j/2}\\
=\exp[O(j(k)\log j(k))]\cdot k^{-j(k)/2}\cdot\exp\biggl(\tfrac{1}{2}\int_0^{j(k)}\!\!\!\! x\log(x/e)\,dx\!\biggr)\\
=\exp\Bigl(-\tfrac{\log^2k}{4\log\log k}+O\bigl(\tfrac{\log^2k}{\log^2(\log k)}\bigr)\Bigr).
\end{multline}

{\bf (b)\/} Suppose that $j\ge  j(k)$. Then, since $j<d$,
\begin{multline*}
\tfrac{k}{j!}\times\bigl|\exp(i(\xi_1/2+j\xi_2/2))-1\bigr|=O\Bigl(\tfrac{k^{1/2+\eps}d}{(j(k)/e)^{j(k)}}\Bigr)\\
=O\bigl(\exp(-\Theta(1/2-\eps)\log k)\bigr).\\
\end{multline*}
Here and below we use the notation $\Theta(1/2-\eps)$ as a shorthand for a quantity which is bounded below by $1/2-\eps$ times an absolute constant.
So, 
%introducing 
%\[
%X_j=\tfrac{k\exp\bigl(-\tfrac{1}{2}\bigl(\tfrac{k-1}{k}\bigr)^2\bigr)}{j!} \bigl(\tfrac{k-1}{k}\bigr)^j,
%\]
we have 
\[
Z_jy_j=X_jy_j+O\bigl(d \|\boldsymbol\xi\|X_jy_j\bigr).
\]
Hence
\begin{multline*}
\sum_{j(k)\le j <d} Z_jy_j=\sum_{j(k)\le j\le d}X_jy_j +O\bigl(d^2\|\boldsymbol\xi\|\tfrac{k}{j(k)!}\bigr)\\
=\sum_{j(k)\le j<d}X_jy_j+ O\bigl(\exp(-\Theta(1/2-\eps)\log k)\bigr),
\end{multline*}
implying that 
\begin{multline*}
\prod_{j(k)\le j< d}\int_0^4\tfrac{\exp\bigl(-Z_j y_j\bigr)}{\sqrt{y_j(4-y_j)}}\,dy_j\\
=\prod_{j(k)\le j<d}\int_0^4\tfrac{\exp\bigl(-X_j y_j\bigr)}{\sqrt{y_j(4-y_j)}}\,dy_j \cdot \bigl[1+ O\bigl(\exp(-\Theta(1/2-\eps)\log k)\bigr)\bigr].
\end{multline*}
Here, using convexity of the exponential function, 
\begin{multline*}
\pi^{-1}\int_0^4\tfrac{\exp\bigl(-X_j y_j\bigr)}{\sqrt{y_j(4-y_j)}}\,dy_j
\ge \exp\biggl(-\tfrac{X_j}{\pi}\int_0^4\sqrt{\tfrac{y}{4-y}}\,dy\biggr)=e^{-2X_j}.\\
\end{multline*}
So,
\begin{multline}\label{3.10}
\prod_{j(k)\le j<d}\pi^{-1}\int_0^4\tfrac{\exp\bigl(-X_j y_j\bigr)}{\sqrt{y_j(4-y_j)}}\,dy_j 
\ge\exp\Bigl(-2\sum_{j\ge j(k)}X_j\Bigr)\\
= \exp\bigl(O(k/j(k)!)\bigr)=k^{O\bigl(\tfrac{\log\log\log k}{\log\log k}\bigr)}.
\end{multline}
%Consequently
%\begin{multline*}
%\prod_{j<d}\int_0^4\tfrac{\exp\bigl(-Z_j y_j\bigr)}{\sqrt{y_j(4-y_j)}}\,dy_j
%=\prod_{j<d}\int_0^4\tfrac{\exp\bigl(-X_j y_j\bigr)}{\sqrt{y_j(4-y_j)}}\,dy_j\cdot [1+O(\Delta_k)],\\
%\end{multline*}
%with $X_j$ defined in \eqref{3.08}. 
Combining \eqref{3.09} and \eqref{3.10} we arrive at
\begin{multline}\label{3.11}
\prod_{j<d}\int_0^4\tfrac{\exp\bigl(-Z_j y_j\bigr)}{\sqrt{y_j(4-y_j)}}\,dy_j 
=\exp\Bigl(-\tfrac{\log^2k}{4\log\log k}+O\bigl(\tfrac{\log^2k}{\log^2(\log k)}\bigr)\Bigr).\\
\end{multline}
Thus,  equations \eqref{3.07} and \eqref{3.11} transform  equation \eqref{3.055} into
\begin{multline}\label{3.12}
\tfrac{1}{z_1^k z_2^{k-1}}\prod_{j< d} H\Bigl(\!\tfrac{z_1z_2^j}{(j!)^2}\Bigr)
=\tfrac{e^{2k}}{x_1^k x_2^{k-1}}\exp\Bigl(-\tfrac{\log^2k}{4\log\log k}+O\bigl(\tfrac{\log^2k}{\log^2(\log k)}\bigr)\Bigr)\\
\times \exp\bigl[-\tfrac{k}{4}\bigl(\xi_1^2+(x_2+x_2^{1/2})\xi_2^2-2\xi_1\xi_2x_2^{1/2}\bigr)\bigr],\\
\end{multline}
uniformly for $\boldsymbol\xi$ with $\|\boldsymbol\xi\|\le k^{-1/2+\eps}$, $\eps\in (0,1/6)$.

By \eqref{3.12}, the contribution to the Cauchy integral on the RHS of \eqref{3.105} is 
\begin{multline*}
\tfrac{e^{2k}}{x_1^k x_2^{k-1}}\exp\Bigl(-\tfrac{\log^2k}{4\log\log k}+O\bigl(\tfrac{\log^2k}{\log^2(\log k)}\bigr)\Bigr)\\
\times\int_{\|\boldsymbol\xi\|\le k^{-1/2+\eps}}\exp\bigl[-\tfrac{k}{4}\bigl(\xi_1^2+(x_2+x_2^{1/2})\xi_2^2-2\xi_1\xi_2x_2^{1/2}\bigr)\bigr]\,d\xi_1 d\xi_2\\
=\tfrac{e^{2k}}{x_1^k x_2^{k-1}}\exp\Bigl(-\tfrac{\log^2k}{4\log\log k}+O\bigl(\tfrac{\log^2k}{\log^2(\log k)}\bigr)\Bigr).\\
\end{multline*}
In combination with \eqref{onestar}, this yields that 
\begin{multline*}
\tfrac{1}{(2\pi i)^2}\!\!\oint\limits_{C_1\times C_2}\!\prod_{j\ge 0}\tfrac{1}{z_1^{k+1} z_2^k}H\Bigl(\tfrac{z_1z_2^j}{(j!)^2}\Bigr)\,dz_1 dz_2\!=\!\tfrac{e^{4k}}{x_1^k x_2^{k-1}}\exp\Bigl(\!-\tfrac{\log^2k}{4\log\log k}+O\bigl(\tfrac{\log^2k}{\log^2(\log k)}\bigr)\!\Bigr)\\
=\tfrac{e^{4k}}{k^{2k}}\exp\Bigl(\!-\tfrac{\log^2k}{4\log\log k}+O\bigl(\tfrac{\log^2k}{\log^2(\log k)}\bigr)\!\Bigr).
\end{multline*}
Therefore, by \eqref{3.105}, $N(k)$ (the total number of pairs of twin trees, each with $k$ vertices, and maximum degree  
$\lfloor2\log k/\log\log k\rfloor$ is given by
\[
N(k)=(k!)^4\cdot\tfrac{e^{4k}}{k^{2k}}\exp\Bigl(\!-\tfrac{\log^2k}{4\log\log k}+O\bigl(\tfrac{\log^2k}{\log^2(\log k)}\bigr)\!\Bigr).
\]
Combining this with \eqref{4}, we obtain that  the expected number of twin trees of size $k=o(n)$ in the random Cayley tree 
with $n$ vertices is at least
\begin{multline*}
\tfrac{1}{n^{n-1}}\binom{n}{k}\binom{n-k}{k}(n-2k)^{n-2k+1}\cdot N(k)\\
=n^2\exp\Bigl(-\tfrac{\log^2k}{4\log\log k}+O\bigl(\tfrac{\log^2k}{\log^2(\log k)}\bigr)\!\Bigr).
\end{multline*}
(Observe that this lower bound depends on $k$ almost like the upper bound \eqref{3.05}, since $\eps''$ can be chosen arbitrarily close to $1/2$ from {\it below\/}.) And the lower bound diverges to infinity if 
\[
k\le\Big\lfloor\! \exp\Bigl((2-\delta)\sqrt{\log n\cdot\log(\log n)}\Bigr)\Big\rfloor.
\]
\end{proof}
{\bf Acknowledgment.\/} We sincerely thank the hard-working referees for helping us to improve the paper, both mathematically and stylistically.

\end{document}